\documentclass[11pt]{article}
\newtheorem{theorem}{Theorem}{\bf}{\it}
\newtheorem{lemma}{Lemma}{\bf}{\it}
\newcommand{\Dom}{\mathop{\mathrm{Dom}}\nolimits}
\newcommand{\sign}{\mathop{\mathrm{sign}}\nolimits}
\newcommand{\im}{\mathop{\mathrm{Im}}\nolimits}
\newcommand{\re}{\mathop{\mathrm{Re}}\nolimits}
\renewcommand{\pmatrix}[1]{\! \left [ \matrix{#1} \right ] \!}
\title{Almost Finite Speed of Propagation for Linear Peridynamics}
\begin{document}
\maketitle
\section{Main Results}
We start with what is essentially a trivial observation about
solutions to the wave equation in one space dimension.
\begin{theorem}\label{t1}
Suppose
\[
  u _ { t t } - c ^ 2 u _ { x x } = 0
\]
and $u$ and $\partial _ t u$ belong to the Schwartz class
$\mathcal S ( \mathbf R )$ when $t = 0$.  For
any $v$ with $| v | \neq c$, any $x _ 0 \in \mathbf R$
and any non-negative integer $l$ there is a constant
$C$ such that
\[
  | \partial _ t u ( t , x _ 0 + v t ) | ^ 2 + c ^ 2 | u _ x ( t , x _ 0 + v t ) | ^ 2
  \le C t ^ { - 2 l }
\]
for all $t > 0$.
\end{theorem}
We will prove the following lemma.
\begin{lemma}
Suppose $J \in L ^ 1 ( \mathbf R )$ is non-negative and even
and that its moments\footnote{All integrals are to be understood as
Lebesgue integrals throughout the paper.}
\[
  \mu _ k = \! \! \int \limits _ { x \in \mathbf R } \! \! x ^ k J ( x ) \, d x
\]
exist for $k \le l + 2$, with $l \ge 0$.
There is a unique $c > 0$ and a unique bounded translation-invariant operator
$D \colon L ^ 2 ( \mathbf R ) \to L ^ 2 ( \mathbf R )$ such that
\begin{enumerate}
\item
\[
  J \star v - \mu _ 0 v = c ^ 2 D ^ 2 v
\]
for all $v$ in $L ^ 2 ( \mathbf R )$.
\item
$- i H D$ is a positive operator, where $H \colon L ^ 2 ( \mathbf R )
\to L ^ 2 ( \mathbf R )$ is the Hilbert transform,
\[
  ( H v ) ( x ) = \lim _ { \epsilon \to 0 ^ + } \! \! \! \! \! \int \limits _ { | x - y | > \epsilon } \! \! \! \! \! \frac { v ( y ) } { x - y } \, d y .
\]
\item
Letting $S _ \lambda \colon L ^ 2 ( \mathbf R )
\to L ^ 2 ( \mathbf R )$ denote scaling by a factor $\lambda$,
\[
  ( S _ \lambda v ) ( x ) = v ( \lambda ^ { - 1 } x ) ,
\]
we have
\[
  \lim _ { h \to 0 ^ + } h ^ { - 1 } S _ h D S _ { h ^ { - 1 } } v
  = v _ x
\]
weakly for any $v \in \Dom ( \partial _ x ) \subset L ^ 2 ( \mathbf R )$.
\end{enumerate}
This $D$ is of the form
\[
  ( D v ) ( x ) = \! \! \int \limits _ { \xi \in \mathbf R } \, \int \limits _ { y \in \mathbf R } \! \! 2 \pi i \psi ( \xi ) e ( \xi ( x - y ) ) v ( y ) \, d \xi \, d y ,
\]
where
\[
  e ( z ) = \exp ( 2 \pi i z ) ,
\]
and $\psi$ is a function uniformly bounded, together with its derivatives up through order~$l$, on~$\mathbf R$.
\end{lemma}
This suggests the following theorem, which we will also
prove.
\begin{theorem}\label{t2}
Suppose that $u$ satisfies the integrodifferential equation\footnote{The convolution of $J$ with $u$ is, of course, to be understood as a convolution in the spatial variable only, for each value of the temporal variable.}
\[
  u _ { t t } + \mu _ 0 u = J \star u
\]
with $J$ as in the lemma.
Suppose that $u$ and $\partial _ t u$ belong to the Schwartz
class $\mathcal S ( \mathbf R )$ when $t = 0$.
For any $v$ with $| v | > c$,
and $x _ 0 \in \mathbf R$ and any non-negative integer $n < N$
there is a constant $C$ such that
\[
  | \partial _ t u ( t , x _ 0 + v t ) | ^ 2 + c ^ 2 | ( D u ) ( t , x _ 0 + v t ) | ^ 2
  \le C t ^ { - 2 l }
\]
for all $t > 0$.  The constant $c$ and operator $D$ are those
of the lemma.
\end{theorem}
Note that this is not an exact analogue of the previous theorem,
because we assume $| v | > c$ rather than $| v | \neq c$.
It is not clear if this difference is essential or if it is
an artifact of the proof below.

The interest in this integrodifferential equation from peridynamics,
a non-local theory of elasticity introduced by Silling~\cite{silling}.
The equation of motion in peridynamics is
\[
  \rho ( x ) u _ { t t } ( t , x ) = \int f ( u ( t , y ) - u ( t , x ) , y - x , x ) \, d y
\]
For homogeneous materials $\rho$ is constant and $f$ is independent
of its third argument. If we assume further that $f$ is linear in
its first argument, which is true approximately for small displacements
in all models and true exactly for small displacements in some models,
then this equation reduces to the one of theorem. One consequence
of the theorem is that small initial displacements remain small for
all time, so we do not leave the domain of validity of the approximation.

Theorem~\ref{t1} is, of course a simple consequence of the explicit
solution formula
\[\begin{array}{r@{\:}c@{\:}l}
  u ( t , x ) & = & \displaystyle \frac 1 2 u ( 0 , x + c t ) + \frac 1 2 u ( 0 , x - c t )
  \cr & & \displaystyle + \frac 1 { 2 c } \! \! \! \! \! \! \! \! \int \limits _ { x - c t < y < x + c t } \! \! \! \! \! \! \! \! \partial _ t u ( 0 , y ) \, d y .
\end{array}\]
The precise details of this formula are unimportant.
It suffices to observe that\footnote{$u$ may be thought of either as a function of two arguments or as a function of its first argument, taking values in functions of the second. It is in this sense that expressions like $u ( t )$ are to be interpreted. Logically the two points of view are equivalent. No ambiguity is possible because one can always count arguments to determine what is meant. The same, of course, applies to derivatives of $u$ as well.}
\[
  u ( t ) = k _ { 0 } ( t ) \star u ( 0 ) + k _ { - 1 } ( t ) \star \partial _ t u ( 0 )
\]
where $k _ 0 ( t ) $ and $k _ { - 1 } ( t )$ are distributions
supported in the interval $[ - c t , c t ]$.
This observation also leads to finite speed of propagation
for the wave equation. It is therefore of interest that
the equation
\[
  u _ { t t } + \mu _ 0 u = J \star u
\]
does not exhibit finite speed of propagation. More precisely,
one has the following theorem.
\begin{theorem}\label{t3}
Suppose $\nu$ belongs to the Schwartz class~$\mathcal S ' ( \mathbf R )$
and that convolution with $\nu$ is positive and self-adjoint in the
sense that
\[
  \left < \nu \star \chi , \chi \right > \ge 0
\]
and
\[
  \left < \nu \star \chi _ 1 , \chi _ 2 \right >
  = \left < \chi _ 1 , \nu \star \chi _ 2 \right > .
\]
for $\chi , \chi _ 1 , \chi _ 2 \in \mathcal S ( \mathbf R )$.
Then there exist families of tempered distribution $k _ { 0 } ( t ),
k _ { - 1 } ( t ) \in \mathcal S ' ( \mathbf R )$ such that the
solution to the initial value problem for the integrodifferential
equation
\[
  u _ { t t } + \nu \star u = 0
\]
has a unique solution, given by
\[
  u ( t ) = k _ { 0 } ( t ) \star u ( 0 ) + k _ { - 1 } ( t ) \star \partial _ t u ( 0 )
\]
For no $t > 0$ are these distributions of compact support
unless $\nu = - ( \alpha \delta '' + i \beta \delta ' + \gamma )$
with either $\alpha = \beta = 0$ and $\gamma \ge 0$
or $\alpha > 0$ and $4 \alpha \gamma + \beta ^ 2 \ge 0$. 
That is the fundamental solution is of compact support only if
the integrodifferential equation above is really a differential
equation
\[
  u _ { t t } - \alpha u _ { x x } - i \beta u _ x - \gamma u = 0 .
\]
\end{theorem}
We should therefore think of Theorem~\ref{t2} as a substitute for
finite speed of propagation, hence the phrase ``almost finite
speed of propagation'' in the title of this paper.

Theorem~\ref{t2} does not provide the sharpest estimates
obtainable for this problem. With more effort it is possible
to obtain the following result.
\begin{theorem}\label{t4}
Suppose that $u$ satisfies
\[
  u _ { t t } + \mu _ 0 u = J \star u ,
\]
where the integrals
\[
  \! \! \int \limits _ { x \in \mathbf R } \! \! x ^ k J ( x ) \, d x
\]
exist for $m \le l + 3$ and the integrals
\[
  \! \!  \int \limits _ { x \in \mathbf R } \! \! x ^ m \partial ^ n _ x u ( 0 , x ) \, d x
\]
and
\[
  \! \! \int \limits _ { x \in \mathbf R } \! \! x ^ m \partial _ t \partial ^ n _ x u ( 0 , x ) \, d x
\]
exist for $m \le l$ and $n \le K + 2$. Then there is a constant $C _ { j , k , l }$ such that
\[
  | \partial ^ j _ t \partial ^ k _ x u ( t , x ) |
  \le C \max ( 1 , 1 + | x | - c | t | ) ^ { - 2 l }
\]
for all $1 \le j \le l$, $k \le K$ and all $t , x \in \mathbf R$.
\[
  | \partial ^ k _ x u ( t , x ) |
  \le C \max ( 1 , 1 + | x | - c | t | ) ^ { - 2 l - 2 }
\]
for all $j \le l$, $k \le K$ and all $t , x \in \mathbf R$.
\end{theorem}
\section{Proof of Theorem~\ref{t3}}
Fourier transforming the integrodifferential equation
\[
  u _ { t t } + \nu u = 0
\]
gives
\[
  F u _ { t t } + \varphi F u = 0 ,
\]
where
\[
  \varphi = F \nu .
\]

By hypothesis,
\[
  \left < \nu \star \chi , \chi \right > \ge 0
\]
and
\[
  \left < \nu \star \chi _ 1 , \chi _ 2 \right >
  = \left < \chi _ 1 , \nu \star \chi _ 2 \right > .
\]
By Plancherel,
\[
  \left < \varphi F \chi , F \chi \right > \ge 0
\]
and
\[
  \left < \varphi \chi _ 1 , \chi _ 2 \right >
  = \left < \chi _ 1 , \varphi \chi _ 2 \right > .
\]
It follows that $\varphi ( \xi )$ is real and non-negative
for all $\xi \in \mathbf R$.

In vector form the differential equation above is
\[
  \partial _ t \pmatrix { F u \cr F \partial _ t u } ( t , \xi )
  = \pmatrix { 0 & 1 \cr - \varphi ( \xi ) } \pmatrix { F u \cr F \partial _ t u } .
\]
Proceeding one frequency at a time,
\[
  \pmatrix { F u \cr F \partial _ t u } ( t , \xi )
  = M ( \sqrt { \varphi ( \xi ) } , t ) \pmatrix { F u \cr F \partial _ t u } ( 0 , \xi )
\]
is the unique solution,
where the matrix valued function $M$ is given by
\[\begin{array}{r@{\:=\:}l}
  M ( \omega , t ) & \displaystyle \exp \left ( \pmatrix { 0 & 1 \cr - \omega ^ 2 & 0 } t \right ) \cr
  & \displaystyle \pmatrix { \cos ( \omega t ) & \omega ^ { - 1 } \sin ( \omega t ) \cr - \omega \sin ( \omega t ) & \cos ( \omega t ) } .
\end{array}\]
Since $M$ is even in $\omega$ the choice of square root is irrelevant.
\[
  \pmatrix { F u \cr F \partial _ t u } ( t , \xi )
  = M ( \sqrt { \varphi ( \xi ) } , t ) \pmatrix { F u \cr F \partial _ t u } ( 0 , \xi ) .
\]
Taking an inverse Fourier transform and a Fourier transform on the other,
we obtain an integral representation for the solution $u$ and its time derivative in terms of the initial data,
\[
  \pmatrix { u \cr \partial _ t u } ( t , x ) = \! \! \! \int \limits _ { \! \! \! \xi \in \mathbf R } \int \limits _ { y \in \mathbf R \! \! \! } \! \! \kappa ( t , x - y , \xi )
\pmatrix { u \cr \partial _ t u } ( 0 , y ) \, d y \, d \xi ,
\]
where
\[
  \kappa ( t , z , \xi ) = e ( \xi z ) M ( \sqrt { \varphi ( \xi ) } , t ) .
\]
It is tempting to reverse the order of integration, obtaining
\[
  \pmatrix { u \cr \partial _ t u } ( t , x ) = \! \! \int \limits _ { y \in \mathbf R } \! \! K ( t , x - y ) \pmatrix { u \cr \partial _ t u } ( 0 , y ) \, d y ,
\]
where
\[
  K ( t , z ) = \! \! \int \limits _ { \xi \in \mathbf R } \! \! \kappa ( t , z , \xi ) \, d \xi ,
\]
but these integrals do not, in general, converge.

To get around this failure of convergence we split off a factor of
\[
  \alpha ( \xi ) = 1 + 4 \pi ^ 2 A \xi ^ 2 
\]
where $A > 0$.
Defining
\[
  L = 1 - A \partial ^ 2 _ x ,
\]
we have
\[
  F L \pmatrix { u \cr \partial _ t u } = \alpha ( \xi ) F \pmatrix { u \cr \partial _ t u }  
\]
Proceeding nearly as before,
\[
  \pmatrix { F u \cr F \partial _ t u } ( t , \xi )
  = \alpha ( \xi ) ^ { - 1 } M ( \sqrt { \varphi ( \xi ) } , t ) \pmatrix { F L u \cr F L \partial _ t u } ( 0 , \xi )
\]
leads to
\[
  \pmatrix { u \cr \partial _ t u } ( t , x ) = \! \! \! \int \limits _ { \! \! \!\xi \in \mathbf R } \int \limits _ { y \in \mathbf R \! \! \! } \! \! \beta ( t , x - y , \xi )
\pmatrix { L u \cr L \partial _ t u } ( 0 , y ) \, d y \, d \xi ,
\]
where
\[
  \beta ( t , z , \xi ) = \alpha ( \xi ) ^ { - 1 } e ( \xi z ) M ( \sqrt { \varphi ( \xi ) } , t ) .
\]
Splitting $e ( \xi z )$ into $\cos ( 2 \pi \xi z )$ and $i \sin ( 2 \pi \xi z )$, we see that the latter is odd and does not contribute to the integral.
Thus we can equally well write
\[
  \pmatrix { u \cr \partial _ t u } ( t , x ) = \! \! \! \int \limits _ { \! \! \!\xi \in \mathbf R } \int \limits _ { y \in \mathbf R \! \! \! } \! \! \lambda ( t , x - y , \xi )
\pmatrix { L u \cr L \partial _ t u } ( 0 , y ) \, d y \, d \xi ,
\]
where
\[
  \lambda ( t , z , \xi ) = \alpha ( \xi ) ^ { - 1 } \cos ( 2 \pi \xi z ) M ( \sqrt { \varphi ( \xi ) } , t ) .
\]
This time we can reverse the order of integration, obtaining
\[
  \pmatrix { u \cr \partial _ t u } ( t , x ) = \! \! \int \limits _ { y \in \mathbf R } \! \! B ( t , x - y  ) \pmatrix { L u \cr L \partial _ t u } ( 0 , y ) \, d y ,
\]
where
\[
  B ( t , z ) = \! \! \int \limits _ { \xi \in \mathbf R } \! \! \lambda ( t , z , \xi ) \, d \xi ,
\]
because the integrand is bounded by a constant times the integrable factor
$\alpha ( \xi ) ^ { - 1 }$, and hence the Lebesgue Dominated Convergence
Theorem applies.
This gives us a representation of the form
\[
  \pmatrix { u \cr \partial _ t u } ( t ) = B ( t ) \star \pmatrix { L u \cr L \partial _ t u } ( 0 ) .
\]
The convolution is to be interpreted in the usual sense, but it
can be equally well considered as the convolution of a tempered
distribution with a test function. Taking this latter interpretation,
we may move the constant coefficient differential operator $L$,
\[
  \pmatrix { u \cr \partial _ t u } ( t ) = K ( t ) \star \pmatrix { u \cr \partial _ t u } ( 0 )  ,
\]
where
\[
  K = L B
\]
provided that we interpret the differentiations and convolutions
in the sense of distributions.

In terms of components,
\[\begin{array}{r@{\:=\:}l}
  u ( t ) & k _ 0 ( t ) \star u ( 0 ) + k _ { - 1 } ( t ) \star \partial _ t u ( 0 ) , \cr
  \partial _ t u ( t ) & k _ 1 ( t ) \star u ( 0 ) + k _ 0 ( t ) \star \partial _ t u ( 0 ) ,
\end{array}\]
\[
  k _ j = L b _ j ,
\]
\[
  b _ j ( t , z ) = \! \! \int \limits _ { \xi \in \mathbf R } \! \! \lambda _ j ( t , z , \xi ) \, d \xi ,
\]
and
\[
  \lambda _ j ( t , z , \xi ) = \alpha ( \xi ) ^ { - 1 } \cos ( 2 \pi \xi z ) (  \varphi ( \xi ) ) ^ { j / 2 } \theta _ j ( t \sqrt { \varphi ( \xi ) } ) 
\]
where
\[
  \theta _ j ( \zeta ) = \cases { ( - 1 ) ^ { j / 2 } \cos \zeta & if $\zeta$ is even , \cr ( - 1 ) ^ { ( j - 1 ) / 2 } \sin \zeta & if $\zeta$ is odd . }
\]
In fact, differentiation under the integral sign shows that we have, for all $j \ge 0$,
\[
  \partial ^ j _ t u ( 0 ) = b _ j ( t ) \star ( L u ) ( 0 ) + b _ { j - 1 } ( t ) \star ( L \partial _ t u ) ( 0 )
\]
and
\[
  \partial ^ j _ t u ( 0 ) = k _ j ( t ) \star u ( 0 ) + k _ { j - 1 } ( t ) \star \partial _ t u ( 0 ) .
\]
The various $b$'s are locally integrable functions, so the convolutions in the first of this pair of equations may be understood in the usual sense, provided the initial data are integrable. The $k$'s are distributions of order at most~2, and the convolutions in the second equation are to be understood in the sense of convolutions.
Even more trivially, we have
\[
  \partial ^ j _ t \partial ^ k _ x u ( t ) = b _ j ( t ) \star ( L \partial ^ k _ x u ) ( 0 ) + b _ { j - 1 } ( t ) \star ( L \partial ^ k _ x \partial _ t u ) ( 0 )
\]
and
\[
  \partial ^ j _ t \partial ^ k _ x u ( t ) = k _ j ( t ) \star \partial ^ k _ x u ( 0 ) + k _ { j - 1 } ( t ) \star \partial ^ k _ x \partial _ t u ( 0 ) .
\]

It remains to determine when $k _ 0 ( t )$ and $k _ { - 1 } ( t )$
are of compact support.
Throughout the remainder of this section $t$ is fixed and positive.

If $k _ j ( t )$ is a distribution of order~$n$ supported in the interval $[ - X , X ]$ then there is a constant $C _ j$ such that
\[
  \left | \left < k _ j , \theta \right > \right | \le C _ j \sum _ { m = 0 } ^ n \max _ { - X \le x \le X } | \theta ^ { ( m ) } ( x ) |
\]
for all smooth functions $\theta$.
There is no loss of generality in assuming $C _ j \ge 1$.
Applying this to the function $\theta ( x ) = e ( - \xi x )$ we see
that $F k _ j ( t )$ is an entire function satisfying.
\[
  | ( F k _ j ( t ) ) ( \xi ) | \le C _ j \sum _ { m = 0 } ^ n ( 2 \pi | \xi | ) ^ m e ^ { 2 \pi X | \im \xi | } .
\]
For $\xi \in \mathbf R$ we have
\[
  ( F k _ 0 ( t ) ) ( \xi ) = \cos ( t \sqrt { \varphi ( \xi ) } ) ,
\]
\[
  ( F k _ { - 1 } ( t ) ) ( \xi ) = \frac { \sin ( t \sqrt { \varphi ( \xi ) } ) } { \sqrt { \varphi ( \xi ) } } ,
\]
and
\[
  \varphi ( \xi ) = \frac { ( F k _ 0 ) ( \xi ) ^ 2 } { 1 - ( F k _ { - 1 } ) ( \xi ) ^ 2 } .
\]
The last of these equations can be used to extend $\varphi$ to a meromorphic function, which must then satisfy the preceding two equations except at poles of $\varphi$.
There are, however, no such poles, because if $\eta \in \mathbf C$ is a pole of $\varphi$ then, by the equation $( F k _ 0 ( t ) ) ( \xi ) = \cos ( t \sqrt { \varphi ( \xi ) } )$, $\sqrt \eta$ is an essential singularity of $F k _ 0 ( t )$, which we already saw was entire.
Thus $\varphi$ is an entire function, but the preceding argument gives no bounds for its size. We will need such bounds on circles of large radius.

If
\[
  | \xi | = R
\]
where $R \ge \frac 1 { \pi }$ then 
\[
  | ( F k _ 0 ( t ) ) ( \xi ) | \le 2 ^ { n + 1 } \pi ^ n C _ j R ^ n e ^ { 2 \pi X R } .
\]
Let $q ( \xi ) = e ( t \sqrt { \varphi ( \xi ) } )$.
There is no reason, at this stage, to believe that the square root can be taken in a continuous manner, so we simply choose a square root arbitrarily at each point.
If $| q | \ge \sqrt 2$ then
\[\begin{array}{r@{\:\le\:}l}
  | q | & 2 ( | q | - | q | ^ { - 1 } ) \le 2 | q + q ^ { - 1 } | 
  \cr & 4 | ( F k _ 0 ( t ) ) | \le 2 ^ { n + 3 } \pi ^ n C _ j R ^ n e ^ { 2 \pi X R } .
\end{array}\]
The same estimate holds trivially if $| q | < \sqrt 2$.
The argument above is equally valid if we replace $q$ by $q ^ { - 1 }$ everywhere, so
\[\begin{array}{r@{\:}c@{\:}l}
  e ^ { 2 \pi | \im \sqrt \varphi ( \xi ) | | t | } & = & \max ( | q | , | q ^ { - 1 } | )
  \cr & \le & 4 | ( F k _ 0 ( t ) ) | \le 2 ^ { n + 3 } \pi ^ n C _ j R ^ n e ^ { 2 \pi X R } ,
\end{array}\]
from which
\[
  | \im \sqrt \varphi ( \xi ) | \le g ( R ) ,
\]
where
\[
  g ( R ) = \frac { 2 X R } { | t | } + \frac n { 2 \pi | t | } \log R + \frac 1 { 2 \pi | t | } \log ( 2 ^ { n + 3 } \pi ^ n C _ j ) .
\]
Now
\[\begin{array}{r@{\:}c@{\:}l}
  \re \varphi ( \xi ) & = & ( \re \sqrt { \varphi ( \xi ) } ) ^ 2 - ( \im \sqrt { \varphi ( \xi ) } ) ^ 2
  \cr & \ge & - ( \im \sqrt { \varphi ( \xi ) } ) ^ 2 \ge - g ( R ) ^ 2 .
\end{array}\]
Define
\[
  \upsilon ( \xi ) = \re \varphi ( \xi ) ,
\]
\[
  \upsilon _ + ( \xi ) = \max ( \upsilon ( \xi ) , 0 ) ,
  \quad \upsilon _ - ( \xi ) = \max ( - \upsilon ( \xi ) , 0 ) ,
\]
so that
\[
  \upsilon ( \xi ) = \upsilon _ + ( \xi ) - \upsilon _ - ( \xi ) ,
  \quad | \upsilon ( \xi ) | = \upsilon _ + ( \xi ) - \upsilon _ - ( \xi ) .
\]
$\upsilon$ is the real part of a holomorphic function, hence harmonic.
By the Mean Value Property
\[
  0 = \upsilon ( 0 ) = \frac 1 { 2 \pi R } \! \int \limits _ { | \xi | = R } \! \upsilon ( \xi ) \, d s ,
\]
from which we obtain
\[
  \frac 1 { 2 \pi R } \! \int \limits _ { | \xi | = R } \! \upsilon _ + ( \xi ) \, d s
  = \frac 1 { 2 \pi R } \! \int \limits _ { | \xi | = R } \! \upsilon _ - ( \xi ) \, d s
\]
and
\[
  \frac 1 { 2 \pi R } \! \int \limits _ { | \xi | = R } \! | \upsilon ( \xi ) | \, d s
  = \frac 2 { 2 \pi R } \! \int \limits _ { | \xi | = R } \! \upsilon _ - ( \xi ) \, d s
\]
Our earlier estimate on $\re \sqrt { \varphi ( \xi ) }$ gives
\[
  \upsilon _ - ( \xi ) \le g ( R ) ^ 2
\]
and hence
\[
  \frac 1 { 2 \pi R } \! \int \limits _ { | \xi | = R } \! | \upsilon ( \xi ) | \, d s \le 2 g ( R ) ^ 2 .
\]
Differentiating the Poisson Formula repeatedly gives
\[
  \partial ^ m _ x \partial ^ n _ y \upsilon ( x + i y ) = \frac 1 { 2 \pi R } \! \! \int \limits _ { | \xi | = R } \! \! P _ { m , n } ( x , y , \re \xi , \im \xi ) \upsilon ( \xi ) \, d s
\]
where
\[
  P _ { m , n } ( x , y , x ' , y ' ) = \frac { p _ { m , n } ( x , y , x ' , y ' ) } { q ( x , y , x ' , y ' ) } ,
\]
$p _ { m , n }$ is a polynomial of degree $m + n + 2$, and
\[
  q ( x , y , x ' , y ' ) = ( x - x ' ) ^ 2 + ( y - y ' ) ^ 2 .
\]
We therefore have
\[
  | \partial ^ m _ x \partial ^ n _ y \upsilon ( x + i y ) |
  \le \frac { 2 g ( R ) ^ 2 \max \limits _ { | \xi | = R } p _ { m , n } ( x , y , \re \xi , \im \xi ) } { ( R ^ 2 - x ^ 2 - y ^ 2 ) ^ { m + n + 1 } }
\]
The numerator is $O ( R ^ { m + n + 4 } )$, so the right hand side
tends to zero as $R$ tends to infinity if $m + n > 2$.
Thus $\upsilon$ is a quadratic polynomial, from which it follows
that $\varphi$ is a quadratic polynomial.
We have already seen that it is real on the real line,
so
\[
  \varphi ( \xi ) = \tilde \alpha \xi ^ 2 + \tilde \beta \xi + \tilde \gamma .
\]
Positivity implies that either $\varphi$ is a non-negative constant
or $\tilde \alpha > 0$ and $4 \tilde \alpha \tilde \gamma - \tilde \beta ^ 2 \ge 0$.
Writing
\[
  \alpha = \frac { \tilde \alpha } { 4 \pi ^ 2 } ,
  \quad \beta = \frac { \tilde \beta } { 2 \pi } ,
  \quad \gamma = - \tilde \gamma ,
\]
we have
\[
  \varphi ( \xi ) = - \alpha ( 2 \pi i \xi ) ^ 2 - i \beta ( 2 \pi i \xi ) - c
\]
or
\[
  \nu = - ( \alpha \delta '' + i \beta \delta ' + \gamma \delta ) .
\]
The conditions on $\tilde a$, $\tilde b$ and $\tilde c$ above are equivalent to
\[
  \alpha > 0 \quad 4 \alpha \gamma + \beta ^ 2 \ge 0 .
\]
\section{Proof of the Lemma}
We start by showing the existence of such $c$ and $D$.
Because $J$ is even any moments of odd order are zero,
including the first moment.
The second moment
\[
  \mu _ 2 = \! \! \int \limits _ { x \in \mathbf R } \! \! x ^ 2 J ( x ) \, d x
\]
exists and is positive.  Let $c$ be the positive solution of
$c$,
\[
  c ^ 2 = \frac { \mu _ 2 } 2 .
\]
The Fourier transform
\[
  ( F J ) ( \xi ) = \! \! \int \limits _ { x \in \mathbf R } \! \! J ( x ) e ( - \xi x ) \, d x ,
\]
is then twice continuously differentiable.
The usual properties of the Fourier transform show that
\[
  \varphi = \mu _ 0 - F J
\]
vanishes at the origin along with its derivative and that it is positive
everywhere else.
Its second derivative is
\[
  \varphi '' ( \xi ) = - ( F J ) '' ( \xi ) = 4 \pi ^ 2 \! \! \int \limits _ { x \in \mathbf R } \! \! x ^ 2 J ( x ) e ( - \xi x ) \, d x .
\]
Using the elementary relation
\[
  e ( - \xi x ) = \cos ( 2 \pi \xi x ) - i \sin ( 2 \pi \xi x )
\]
and the evenness of $J$,
\[
  \varphi '' ( \xi ) = 4 \pi ^ 2 \! \! \int \limits _ { x \in \mathbf R } \! \! x ^ 2 J ( x ) \cos ( 2 \pi \xi x ) \, d x
\]
In particular,
\[
  \varphi '' ( 0 ) = 4 \pi ^ 2 \mu _ 2 > 0 .
\]
Since $| \cos ( 2 \pi \xi x ) | \le 1$ everywhere we have
\[
  | \varphi '' ( \xi ) | \le 4 \pi ^ 2 \mu _ 2
\]
for all $\xi \in \mathbf R$.

We define $\psi$ by
\[
  \psi ( \xi ) ^ 2 = 2 \varphi ( \xi ) / \varphi '' ( 0 ) ,
\]
taking the positive square root for positive values of~$\xi$ and the
negative square root for negative values. This choice clearly makes
$\psi$ an odd function.

A simple calculation shows that
\[
  \psi ' ( 0 ) = 1
\]
and
\[
  \psi ' ( \xi ) ^ 2 = \frac { \varphi ' ( \xi ) ^ 2 } { 2 \varphi '' ( 0 ) \varphi ( \xi ) }
\]
for $\xi \neq 0$.
From $\varphi '' ( \xi ) \le \varphi '' ( 0 )$ it follows that
\[
  \varphi ( \eta ) \le \varphi ( \xi ) + \varphi ' ( \xi ) ( \eta - \xi )
    + \frac 1 2 \varphi '' ( 0 ) ( \eta - \xi ) ^ 2
\]
for all $\xi , \eta \in \mathbf R$.  Since $\varphi ( \eta ) \ge 0$,
\[
  \varphi ( \xi ) + \varphi ' ( \xi ) ( \eta - \xi )
    + \frac 1 2 \varphi '' ( 0 ) ( \eta - \xi ) ^ 2 \ge 0 .
\]
Taking
\[
  \eta = \xi + \frac { \varphi ' ( \xi ) } { \varphi '' ( 0 ) }
\]
gives
\[
  \varphi ( \xi ) - \frac 1 2 \frac { \varphi ' ( \xi ) ^ 2 } { \varphi '' ( 0 ) } \ge 0
\]
or
\[
  \psi ' ( \xi ) ^ 2 \le 1 .
\]

In fact, all the derivatives of $\psi$ through order~$l$ are bounded.
By hypothesis $J$ has moments of order up to $l + 2$.
In other words, the integrals
\[
  \mu _ k = \! \! \int \limits _ { x \in \mathbf R } \! \! x ^ m J ( x ) \, d x
\]
exist for $k \le l + 2$. This is equivalent to the existence of the integrals
\[
  \nu _ k \! \! \int \limits _ { x \in \mathbf R } \! \! | x ^ m J ( x ) | \, d x .
\]
Since $J$ is non-negative $\nu _ k = \mu _ k$ if $k$ is even, but if $k$ is odd then $\nu _ k$ is positive while $\mu _ k$ is zero.
In any case
\[
  \varphi ^ { ( k ) } ( \xi ) = - ( 2 \pi i ) ^ k \! \! \int \limits _ { x \in \mathbf R } x ^ k J ( x ) e ( - \xi x ) \, d x
\]
for $1 \le k \le l + 2$, so
\[
  | \varphi ^ { ( k ) } ( \xi ) | \le ( 2 \pi ) ^ k \nu _ k .
\]
Now
\[
  \lim _ { \xi \to \pm \infty } = \mu _ 0
\]
by the Riemann-Lebesgue Lemma. There is therefore a $\Xi$ such that
\[
  | \varphi ( \xi ) - \mu _ 0 | \le \frac { \mu _ 0 } 2
\]
for $| \xi | \ge \Xi$.
\[
  \psi ( \xi ) = \pm \sqrt { \frac { 2 \mu _ 0 } { \varphi '' ( 0 ) } }
    \sqrt { \frac { \varphi ( \xi ) } { \mu _ 0 } } .
\]
The first factor is constant. The second is the composition of the
square root function, which is smooth on the compact set $[ 1 / 2 , 3 / 2 ]$,
with $\varphi$, which is uniformly bounded along with its derivatives
up through order~$l$  everywhere
and takes values in $( 1 / 2 , 3 / 2 )$ for $| \xi | \ge \Xi$.
Thus $\psi$ is uniformly bounded along with its derivatives through
order~$l + 2$ for $| \xi | \ge \Xi$.

Near zero we have to proceed differently. Integration by
parts twice yields the relation
\[
  \int _ 0 ^ \xi ( \xi - \eta ) \varphi '' ( \eta ) \, d \eta = \varphi ( \xi ) .
\]
Making the change of variable $\eta = \tau \xi$,
\[
  \varphi ( \xi ) = \xi ^ 2 \int _ 0 ^ 1 ( 1 - \tau ) \varphi '' ( \tau \xi ) \, d \tau ,
\]
from which we obtain
\[
  \psi ( \xi ) = \xi \sqrt { 2 \int _ 0 ^ 1 ( 1 - \tau ) \frac { \varphi '' ( \tau \xi ) } { \varphi '' ( 0 ) } \, d \tau } .
\]
Since $\varphi ''$ is continuous we can find an $\epsilon > 0$ such that
\[
  | \varphi '' ( \eta ) - \varphi '' ( 0 ) | < \frac 1 2 | \varphi '' ( 0 ) |
\]
for all $| \xi | < \epsilon$.
For such $\xi$ we have
\[
  \frac 1 2 < 2 \int _ 0 ^ 1 ( 1 - \tau ) \frac { \varphi '' ( \tau \xi ) } { \varphi '' ( 0 ) } \, d \tau < \frac 3 2 .
\]
Repeated differentiation under the integral sign shows that this quantity
and its derivatives up through order~$l$ are uniformly bounded.
As before, this is true of the square root as well.
Multiplication by $\xi$ does not change this, since $\xi$ and all
its derivatives are uniformly bounded in $[ - \epsilon , \epsilon ]$.

Finally, we are left with the intervals $[ \epsilon , \Xi ]$ and
$[ - \Xi , - \epsilon ]$.
The $k$'th derivative is bounded by induction on~$k$.
Differentiating
\[
  \xi ^ 3 \int _ 0 ^ 1 \frac { \varphi ''' ( \tau \xi ) } 2 ( 1 - \tau ) ^ 2 \,  d \tau
\]
$k$ times gives
\[
  \psi ( \xi ) \psi ^ { ( k ) } ( \xi ) = \frac { \varphi ^ { ( k ) } ( \xi ) } { \varphi '' ( 0 ) } - \frac 1 2 \sum _ { j = 1 } ^ { k - 1 } \psi ^ { ( k ) } ( \xi ) \psi ^ { ( n - j ) } ( \xi )
\]
As long as $k \le l + 1$ the right hand side is uniformly bounded by the
induction hypothesis. $\psi$ is non-zero in these intervals and, because
the intervals are compact, is bounded away from zero. Therefore
$\psi ^ { ( k ) }$ is uniformly bounded, thus recovering our inductive
hypothesis.

We define the operator $D$ to be multiplication of the
Fourier transform by $\psi$, multiplied by a factor of $2 \pi i$ to make
the result real,
\[
  ( F D v ) ( \xi ) = 2 \pi i \psi ( \xi ) ( F v ) ( \xi )
\]
or
\[
  ( D v ) ( x ) = \! \! \int \limits _ { \xi \in \mathbf R } \, \int \limits _ { y \in \mathbf R } \! \! 2 \pi i \psi ( \xi ) e ( \xi ( x - y ) ) v ( y ) \, d \xi \, d y .
\]
Note that since $D$ is a Fourier multiplier it is automatically
translation invariant.
It is tempting reverse the order of the integrals and write
\[
  D v = Q \star v
\]
where
\[
  Q ( x ) = 2 \pi i \! \! \int \limits _ { \xi \in \mathbf R } \! \! \psi ( \xi ) e ( \xi x ) \, d x .
\]
This integral, however, does not converge.  It can be
given a meaning in the sense of distributions.  The resulting
distribution $Q$ is regular and integrable away from zero, but
has a singularity like $1 / x$ there.  Convolution of $Q$ can
therefore be given a meaning through principal value integrals,
as is done with the Hilbert transform.  This, however, does not
seem to be worth the effort.

Next we note that
\[
  ( F S _ { h ^ { - 1 } } v ) ( \xi ) = h ^ { - 1 } ( F v ) ( h ^ { - 1 } \xi ) ,
\]
\[
  ( F D S _ { h ^ { - 1 } } v ) ( \xi ) = 2 \pi i h ^ { - 1 } \psi ( \xi ) ( F v ) ( h ^ { - 1 } \xi ) ,
\]
\[
  ( F S _ h D S _ { h ^ { - 1 } } v ) ( \xi ) = 2 \pi i \psi ( h \xi ) ( F v ) ( \xi ) ,
\]
and
\[
  ( F S _ h D S _ { h ^ { - 1 } } v ) ( \xi ) = 2 \pi i h ^ { - 1 } \psi ( h \xi ) ( F v ) ( \xi ) .
\]
Because $\psi ( 0 ) = 0$ and $\psi ' ( 0 )$ we have
\[
  \lim _ { h \to 0 ^ + } h ^ { - 1 } \psi ( h \xi ) = \xi .
\]
From $| \psi ' ( \xi ) \le 1$ it follows that
\[
  | h ^ { - 1 } \psi ( h \xi ) | \le | \xi | .
\]
If then $w \in L ^ 2 ( \mathbf R )$ and $v \in \Dom ( \partial _ x )$ then
\[
  \lim _ { h \to 0 ^ + } \! \! \int \limits _ { \xi \in \mathbf R } \! \! h ^ { - 1 } ( F S _ h D S _ { h ^ { - 1 } } v ) ( \xi ) \overline { ( F w ) ( \xi ) } \, d \xi
\]
is
\[
  \int \limits _ { \xi \in \mathbf R } \! \! 2 \pi i \xi ( F v ) ( \xi ) \overline { ( F w ) ( \xi ) } \, d \xi ,
\]
by Lebesgue Dominated Convergence.
Since the Fourier transform is an isometry this is equivalent to the
statement that
\[
  \lim _ { h \to 0 ^ + } \left < h ^ { - 1 } S _ h D S _ { h ^ { - 1 } } v
  , w \right >
  = \left < v _ x , w \right > ,
\]
which is the weak convergence promised in the lemma.

Finally we observe that
\[
  - i ( F H D v ) ( \xi ) = \pi \sign \xi \, \psi ( \xi ) ( F V ) ( \xi ) .
\]
The Fourier multiplier $\pi \sign \xi \, \psi ( \xi )$ is positive
almost everywhere, so $- i H D$ is a positive operator.
This concludes the proof of the existence of $c$ and $D$.

\section{Proof of Theorem~\ref{t2}}
To prove this, first note that, by the usual formula for the Fourier
transform of a convolution,
\[
  \mu _ 0 u - J \star u = - c ^ 2 D ^ 2 u .
\]
Our integrodifferential equation can therefore be rewritten
as
\[
  u _ { t t } - c ^ 2 D ^ 2 u = 0 .
\]
Introducing the quantities
\[
  p = \partial _ t u , \quad q = c D u ,
\]
our second order equation is equivalent to the first order
system
\[
  p _ t = c D q , \quad q _ t = c D p .
\]
Defining,
\[
  w ^ \pm = p \pm q ,
\]
we find
\[
  w ^ \pm _ t = \pm c D w ^ \pm .
\]
Taking Fourier transforms,
\[
  F w ^ \pm _ t = \pm 2 \pi i c \psi F w ^ \pm ,
\]
with solution
\[
  ( F w ^ \pm ) ( t , \xi ) = e ( \pm c \psi ( \xi ) t ) ( F w ^ \pm ) ( 0 , \xi )
\]
Taking inverse Fourier transforms,
\[
  w ^ \pm ( t , x ) = \! \! \int \limits _ { \xi \in \mathbf R } \! \! e ( \xi x \pm c \psi ( \xi ) t ) ( F w ^ \pm ) ( 0 , \xi ) \, d \xi .
\]
In particular,
\[
  w ^ \pm ( t , x _ 0 + v t ) = \! \! \int \limits _ { \xi \in \mathbf R }  \! \! e ( ( \xi x _ 0 + s ^ \pm t ) ( F w ^ \pm ) ( 0 , \xi ) \, d \xi .
\]
where
\[
  s ^ \pm ( \xi ) = v \xi \pm c \psi ( \xi ) .
\]
The integrodifferential equation is symmetric under reflections,
so we may assume without loss of generality that $v > c$.
Then
\[
  \partial _ \xi s ^ \pm \ge v - c > 0 .
\]
We now integrate by parts repeatedly,
\[
  w ^ \pm ( t , x _ 0 + v t )
  = t ^ { - l } \! \! \int \limits _ { \xi \in \mathbf R } \! \! r _ n ( t , \xi ) ) e ( s ^ \pm ( \xi ) t ) \, d \xi ,
\]
where $r _ h$ is given inductively by
\[
  r _ 0 ( t , \xi ) = e ( \xi x _ 0 ) ( F w ^ \pm ) ( 0 , \xi ) ,
\]
\[
  r _ { h + 1 } ( t , \xi ) = - \partial _ \xi ( ( 2 \pi i \partial _ \xi s ^ \pm ( \xi ) ) ^ { - 1 }  r _ h ( t , \xi ) )
\]
Now, by induction on $h$,
\[
  \| r _ h ^ \pm \| _ { L ^ 1 ( \mathbf R ) } < \infty .
\]
Also $w ^ \pm ( 0 , x )$ belongs to $\mathcal S ( \mathbf R )$, so
$( F w ^ \pm ) ( 0 , \xi )$ is also in $\mathcal S ( \mathbf R )$.
It follows that
\[
  | w ^ \pm ( t , x _ 0 + v t ) | \le C ^ \pm t ^ { - 2 l } .
\]
The parallelogram identity
\[
  p ^ 2 + q ^ 2 = \frac 1 2 [ ( w ^ + ) ^ 2 + ( w ^ - ) ^ 2 ]
\]
then gives the required estimate.
\section{Proof of Theorem \ref{t4}}
Starting from the equations
\[
  b _ j ( t , z ) = \! \! \int \limits _ { \xi \in \mathbf R } \! \! \lambda _ j ( t , z , \xi ) \, d \xi
\]
and
\[
  \lambda _ j ( t , z , \xi ) = \alpha ( \xi ) ^ { - 1 } \cos ( 2 \pi \xi z ) ( 2 \pi c \psi ( \xi ) ) ^ j \theta _ j ( 2 \pi c t \psi ( \xi ) ) 
\]
of the previous section
and using the trigonometric identity
\[
  \theta _ j ( x ) \cos ( y ) = \frac 1 2 \left ( \theta _ j ( x + y ) + \theta _ j ( x + y ) \right )
\]
to write
\[
  b _ j = b ^ + _ j + b ^ - _ j
\]
where
\[
  b ^ \pm _ j ( t , z ) = \frac 1 2 \! \! \int \limits _ { \xi \in \mathbf R } \! \! \alpha ( \xi ) ^ { - 1 } ( 2 \pi \psi ( \xi ) c t ) ^ j \theta _ j ( 2 \pi \sigma ^ \pm ( t , z , \xi ) ) \, d \xi
\]
and
\[
  \sigma ^ \pm ( t , z , \xi ) = \psi ( \xi ) c t \pm \xi z
\]
In this formula we integrate by parts $N$ times, using the relation
\[
  \partial _ \zeta ^ N \theta _ { j - N } ( \zeta ) = \theta _ j ( \zeta ) .
\]
It is easy to see that the boundary terms all vanish.
We thus find
\[
  b ^ \pm _ j ( t , z ) = \! \! \int \limits _ { \xi \in \mathbf R } \! \! q ^ \pm _ { j , N } ( t , z , \xi ) \theta _ { j - N } ( 2 \pi \sigma ^ \pm ( t , z , \xi ) ) \, d \xi
\]
where the $q$'s are given inductively by
\[
  q ^ \pm _ { j , 0 } ( t , z , \xi ) = \frac 1 2 \alpha ( \xi ) ^ { - 1 } ( 2 \pi c \psi ( \xi ) ) ^ j ,
\]
\[
  q ^ \pm _ { j , n + 1 } ( t , z , \xi ) = - \partial _ \xi \left [ \rho ^ \pm ( t , z , \xi ) ^ { - 1 } q ^ \pm _ { j , n } ( t , z , \xi ) \right ] ,
\]
where
\[
  \rho ^ \pm ( t , z , \xi ) = 2 \pi \partial _ \xi \sigma ^ \pm ( t , z , \xi ) = 2 \pi ( c t \psi ' ( \xi ) \pm z ) .
\]
If we define functions $\gamma ^ \pm _ { j , n , h }$ for $0 \le h \le n$ by
\[
  \gamma ^ \pm _ { j , 0 , 0 } ( \xi ) = \frac 1 2 \alpha ( \xi ) ^ { - 1 } ( 2 \pi c \psi ( \xi ) ) ^ j
\]
and\footnote{The terms on the right where the last subscript is out of range,
\textit{i.e.} the first term for $k = 0$ and the second for $k = n + 1$, are to be interpreted as zero.}
\[
  \gamma ^ \pm _ { j , n + 1 , h } ( \xi )
  = 2 ( n + h ) \pi c t \psi '' ( \xi ) \gamma ^ \pm _ { j , n , h - 1 } ( \xi ) - \partial _ \xi \gamma ^ \pm _ { j , n , k } ( \xi )
\]
then a simple induction on $n$ shows that
\[
  q ^ \pm _ { j , n } ( t , z , \xi ) = \sum _ { h = 0 } ^ n \gamma ^ \pm _ { j , n , h } ( \xi ) \rho ^ \pm ( t , z , \xi ) ^ { - n - h } .
\]
By induction, $\gamma ^ \pm _ { j , n , k }$ is a linear combination
of terms consisting of derivatives of $\alpha$ multiplied by
products of derivatives of $\psi$. These latter derivatives are
of order at most $n + 1$.
Thus, if $n \le l - 1$,
\[
  C _ { j , n , h } =
  2 \| \gamma ^ \pm _ { j , n , h } \| _ { L ^ 1 ( \mathbf R ) } < \infty .
\]
Since
\[
  | \theta _ { j - N } ( \xi ) | \le 1
\]
and
\[
  | \rho ^ \pm ( t , z , \xi ) | \ge | z | - c | t |
\]
for all $\xi$ we conclude that
\[
  | b ^ \pm _ j ( t , z ) | \le \frac 1 2 \sum _ { h = 0 } ^ n C _ { j , n , h } ( | z | - c | t | ) ^ { - n - h }
\]
if $| z | > c | t |$ and hence
\[
  | b _ j ( t , z ) | \le \sum _ { h = 0 } ^ n C _ { j , n , h } ( | z | - c | t | ) ^ { - n - h } .
\]

Special care is required for $j = - 1$. The prove given above fails for
$j < 0$ because $\psi ( \xi ) ^ { - j }$ and its derivatives are no
longer bounded. $j < -1$ is irrelevant, but we need $j = - 1$ for
the proof of Theorem~\ref{t4}. Fortunately
\[
  \partial _ t b _ { - 1 } ( t , z ) = b _ 0 ( t , z )
\]
and
\[
  b _ { - 1 } ( 0 , z ) = 0
\]
so we can simply integrate the estimates for $j = 0$ to obtain
\[
  | b _ { - 1 } ( t , z ) | \le \sum _ { h = 0 } ^ n \frac 1 { ( n + h - 1 ) c } C _ { 0 , n , h } ( | z | - c | t | ) ^ { - n - h + 1 } .
\]

This estimate and the ones obtained previously are, of course, only
useful when $| z | - c | t |$ is positive and reasonably large.
For $| z | - c | t |$ small or negative it is best to return
to the relation
\[
  b _ j ( t , z ) = \! \! \int \limits _ { \xi \in \mathbf R } \! \! \lambda _ j ( t , z , \xi )
\]
and note that the integrand $\lambda _ j ( t , z , \xi )$
consists of the integrable factor
$\alpha ( \xi ) ^ { - 1 }$ multiplied by something uniformly
bounded, so we have the trivial estimate
\[
  | b _ j ( t , z ) | \le C _ j .
\]
We can combine this with the estimates obtained previously to
find
\[
  | b _ j ( t , z ) | \le C _ { j , m } \min ( 1 , | z | - c | t | ) ^ { - m }
\]
for all $t$ and $z$, where $m = n - 1$ if $j = - 1$ and $m = n$ if $j \ge 0$.

Now we use the relation
\[
  \partial ^ j _ t \partial ^ k _ x u ( t ) = b _ j ( t ) \star ( L \partial ^ k _ x u ) ( 0 ) + b _ { j - 1 } ( t ) \star ( L \partial ^ k _ x \partial _ t u ) ( 0 )
\]
derived in the previous section.
Writing this as
\[
  \partial ^ j _ t \partial ^ k _ x u ( t , x ) = U _ { j , k } ( t , x ) + U ' _ { j , k } ( t , x ) ,
\]
where
\[
  U _ { j , k } ( t , x ) = \! \! \int \limits _ { y \in \mathbf R } \! \! b _ j ( t , x - y ) ( L \partial ^ k _ x u ) ( 0 , y ) \, d y
\]
and
\[
  U ' _ { j , k } ( t , x ) = \! \! \int \limits _ { y \in \mathbf R } \! \! b _ { j - 1 } ( t , x - y ) ( L \partial ^ k _ x u ) ( 0 , y ) \, d y ,
\]
we have the trivial estimates
\[
  | U _ { j , k } ( t , x ) | \le W _ { k , l } Z _ { j , l }
\]
and
\[
  | U ' _ { j , k } ( t , x ) | \le W ' _ { k , l } Z _ { j - 1 , l }
\]
where
\[
  W _ { k , l } = \! \! \int \limits _ { y \in \mathbf R } \! \! ( 1 + | y | ) ^ l | \partial ^ k _ x u ( 0 , y ) | \, d y
\]
\[
  W ' _ { k , l } = \! \! \int \limits _ { y \in \mathbf R } \! \! ( 1 + | y | ) ^ l | \partial ^ k _ x \partial _ t u ( 0 , y ) | \, d y ,
\]
and
\[
  Z _ { j , l } = \max _ { y \in \mathbf R } | b _ { j - 1 } ( t , x - y ) | ( 1 + | y | ) ^ { - l }
\]
Now
\[
  W _ { k , l } \le \sum _ { m = 0 } ^ l \frac { l ! } { m ! ( l - m ) ! }
    ( V _ { k , m } + A V _ { k + 2 , m } )
\]
and
\[
  W ' _ { k , l } \le \sum _ { m = 0 } ^ l \frac { l ! } { m ! ( l - m ) ! }
    ( V ' _ { k , m } + A V ' _ { k + 2 , m } ) ,
\]
where
\[
  V _ { k , m } = \! \! \int \limits _ { y \in \mathbf R } \! \! | y | ^ m | \partial ^ k _ x u ( 0 , y ) | \, d y
\]
and
\[
  V ' _ { k , m } = \! \! \int \limits _ { y \in \mathbf R } \! \! | y | ^ m | \partial ^ k _ x \partial _ t u ( 0 , y ) | \, d y .
\]
By hypothesis $V _ { k , m }$ and $V ' _ { k , m }$ are finite for $k \le K + 2$ and $m \le l$, so $W _ { j , l }$ and $W ' _ { j , l }$ are finite.

Also
\[
  Z _ { j , l } \le C _ { j , l } \max _ { y \in \mathbf R } \min ( 1 , | x - y | - c | t | ) ^ { - l } ( 1 + | y | ) ^ { - l }
\]
or, evaluating the right hand side,
\[
  Z _ { j , l } \le \cases {
    C _ { j , l } & if $| x | \le c | t | + 1$, \cr
    2 ^ { 2 l } C _ { j , l } ( 1 + | x | - c | t | ) ^ { - 2 l }
      & if $| x | \ge c | t | + 1$. }
\]
From this it follows that
\[
  Z _ { j , l } \le 2 ^ { 2 l } C _ { j , l } \min ( 1 , | x | - c | t | ) ^ { - 2 l }
\]
for all $t , x \in \mathbf R$ and hence that
\[
  | U _ { j , k } ( t , x ) \le 2 ^ { 2 l } C _ { j , l } W _ { k , l } \min ( 1 , | x | - c | t | ) ^ { - 2 l }
\]
The same is true for $U ' _ { j , k } ( t , x )$, with $W ' _ { k , l }$ in
place of $W _ { k , l }$ and $C _ { j - 1 , l }$ in place of $C _ { j , l }$,
except possibly for $j = 0$ where we have to
substitute our weaker estimate for $b _ { - 1 }$, obtaining
\[
  | U ' _ { 0 , k } ( t , x ) \le 2 ^ { 2 l - 2 } C _ { - 1 , l } W ' _ { k , l } \min ( 1 , | x | - c | t | ) ^ { - 2 l - 2 } .
\]
Combining all these estimates,
\[
  | \partial ^ j _ t \partial ^ k _ x u ( t , x ) | \le C _ { j , k , l } \min ( 1 , | x | - c | t | ) ^ { - 2 l }
\]
if $0 < j \le l$ and
\[
  | \partial ^ k _ x u ( t , x ) | \le C _ { 0 , k , l } \min ( 1 , | x | - c | t | ) ^ { - 2 l - 2 } .
\]
This concludes the proof of Theorem~\ref{t4}.
\bibliography{afsp}
\bibliographystyle{siam}
\end{document}